\documentclass{article}
\usepackage{amsfonts}
\usepackage{amsmath}

\newenvironment{proof}[1][Proof]{\noindent\textbf{#1.} }{\ \rule{0.5em}{0.5em}}
\input{tcilatex}

\begin{document}

\section{Introduction}

\qquad It is well known(see [CH], [K], or [MP]) that there is a unique
ellipse inscribed in a given triangle, $T$, tangent to the sides of $T$ at
their respective midpoints. This is often called the midpoint or Steiner
inellipse, and it can be characterized as the inscribed ellipse having
maximum area. In addition, one has the following inequality. If $E$ denotes
any ellipse inscribed in $T$, then 
\begin{equation}
\dfrac{\text{Area}\left( E\right) }{\text{Area}\left( T\right) }\leq \dfrac{%
\pi }{3\sqrt{3}},  \tag{trineq}
\end{equation}%
with equality if and only if $E$ is the midpoint ellipse. In [MP] the
authors also discuss a connection between the Steiner ellipse and the
orthogonal least squares line for the vertices of $T$.

\textbf{Definition:} The line, $\pounds $, is called a line of best fit for $%
n$ given points $z_{1},...,z_{n}$ in $%
\mathbb{C}
$, if $\pounds $ minimizes $\tsum\limits_{j=1}^{n}d^{2}\left( z_{j},l\right) 
$ among all lines $l$ in the plane. Here $d\left( z_{j},l\right) $ denotes
the distance(Euclidean) from $z_{j}$ to $l$.

\qquad In [MP] the authors prove the following results, some of which are
also proven in

\textbf{Theorem A: }Suppose $z_{j}$ are points in $%
\mathbb{C}
,g=\dfrac{1}{n}\tsum\limits_{j=1}^{n}z_{j}$ is the centroid, and $%
Z=\tsum\limits_{j=1}^{n}\left( z_{j}-g\right)
^{2}=\tsum\limits_{j=1}^{n}z_{j}^{2}-ng^{2}$.

(a) If $Z=0$, then every line through $g$ is a line of best fit for the
points $z_{1},...,z_{n}$.

(b) \ If $Z\neq 0$, then the line, $\pounds $, thru $g$ that is parallel to
the vector from $(0,0)$ to $\sqrt{Z}$ is the unique line of best fit for $%
z_{1},...,z_{n}$.

\textbf{Theorem B: }If $z_{1},z_{2},$ and $z_{3}$ are the vertices of a
nonequilateral triangle, $T$, and $\pounds $ minimizes $\tsum%
\limits_{k=1}^{3}d^{2}\left( z_{k},l\right) $ among all lines $l$ in the
plane, then the major axis of the Steiner ellipse lies on $\pounds $.

A proof of Theorem B goes all the way back to Coolidge in 1913(see [MP]).

The purpose of this paper is to attempt to generalize (trineq) and Theorem B
to ellipses inscribed in convex quadrilaterals. Many of the results in this
paper use results from two earlier papers of the author about ellipses
inscribed in quadrilaterals. In particular, in [H1] we proved the following
results.

\textbf{Theorem C: }Let D\ be a convex quadrilateral in the $xy$ plane and
let $M_{1}$ and $M_{2}$ be the midpoints of the diagonals of D. Let $Z$ be
the open line segment connecting $M_{1}$ and $M_{2}$. If $(h,k)\in Z$ then
there is a unique ellipse with center $(h,k)$ inscribed in D.

\textbf{Theorem D:} Let D\ be a convex quadrilateral in the $xy$ plane. Then
there is a unique ellipse of maximal area inscribed in D.

\qquad The following general result about ellipses is essentially what
appears in [W], except that the cases with $A=B$ were added by the author.

\textbf{Lemma E: }Let $E$ be an ellipse with equation $%
Ax^{2}+By^{2}+2Cxy+Dx+Ey+F=0$, and let $\phi $\ denote the counterclockwise
angle of rotation from the line thru the center parallel to the $x$\ axis to
the major axis of $E$. Then

\begin{equation*}
\phi =\left\{ 
\begin{array}{ll}
0 & \text{ if }C=0\text{ and }A<B \\ 
\dfrac{\pi }{2} & \text{ if }C=0\text{ and }A>B \\ 
\dfrac{1}{2}\cot ^{-1}\left( \dfrac{A-B}{2C}\right) & \text{if }C\neq 0\ 
\text{and }A<B \\ 
\dfrac{\pi }{2}+\dfrac{1}{2}\cot ^{-1}\left( \dfrac{A-B}{2C}\right) & \text{%
if }C\neq 0\ \text{and }A>B \\ 
\dfrac{\pi }{4} & \text{if }C<0\ \text{and }A=B \\ 
\dfrac{3\pi }{4} & \text{if }C>0\ \text{and }A=B%
\end{array}%
\right.
\end{equation*}

\qquad In [H2], we also derived the following results about ellipses
inscribed in parallelograms.

\textbf{Lemma F:} Let $Z$\ be the rectangle with vertices $%
(0,0),(l,0),(0,k), $and $(l,k),$ where $l,k>0$.

(A) The general equation of an ellipse, $E$, inscribed in $Z$\ is given by 
\begin{equation*}
k^{2}x^{2}+l^{2}y^{2}-2l\left( k-2v\right)
xy-2lkvx-2l^{2}vy+l^{2}v^{2}=0,0<v<k.
\end{equation*}

(B) The corresponding points of tangency of $\Psi $ are 
\begin{equation*}
\left( \dfrac{lv}{k},0\right) ,(0,v),\left( \dfrac{l}{k}\left( k-v\right)
,k\right) ,\text{and\ }(l,k-v).
\end{equation*}

(C) If $a$ and $b$ denote the lengths of the semi--major and semi--minor
axes, respectively, of $E$, then 
\begin{eqnarray*}
a^{2} &=&\dfrac{2l^{2}\left( k-v\right) v}{k^{2}+l^{2}-\sqrt{%
(k^{2}+l^{2})^{2}-16l^{2}\left( k-v\right) v}}\text{ and} \\
b^{2} &=&\dfrac{2l^{2}\left( k-v\right) v}{k^{2}+l^{2}+\sqrt{%
(k^{2}+l^{2})^{2}-16l^{2}\left( k-v\right) v}}.
\end{eqnarray*}

\textbf{Proposition G:} Let P\ be the parallelogram with vertices $%
O=(0,0),P=(l,0),Q=(d,k),$ and $R=(d+l,k)$, where $l,k>0,d\geq 0$. The
general equation of an ellipse, $E$, \textit{inscribed} in P\ is given by%
\begin{gather*}
k^{3}x^{2}+\left( k(d+l)^{2}-4dlv\right) y^{2}-2k\left( k(d+l)-2lv\right) xy
\\
-2k^{2}lvx+2klv\left( d-l\right) y+kl^{2}v^{2}=0,0<v<k.
\end{gather*}

Of course, any two triangles are affine equivalent, while the same is not
true of quadrilaterals--thus it is not surprising that not all of the
results about ellipse inscribed in triangles extend nicely to
quadrilaterals. For example, there is not necessarily an ellipse inscribed
in a given quadrilateral, D, which is tangent to the sides of D\ at their
respective midpoints. There is such an ellipse, which we call the midpoint
ellipse, when D\ is a parallelogram(see Proposition 2.2). We are able to
prove an inequality(see Theorem 2.1), similar to (trineq), which holds for
all convex quadrilaterals. If $E$ is any ellipse inscribed in a
quadrilateral, D, then $\dfrac{\text{Area}\left( E\right) }{\text{Area(D)}}%
\leq \dfrac{\pi }{4}$, and equality holds if and only if D\ is a
parallelogram and $E$ is the midpoint ellipse.

Not suprisingly, Theorem B also does not extend in general to ellipses
inscribed in convex quadrilaterals. However, such a characterization does
hold again when D\ is a parallelogram. We prove in Theorem 3.1 that the foci
of the unique ellipse of maximal area inscribed in a parallelogram, D, lie
on the orthogonal least squares line for the vertices of D. It is also well
known that if $p(z)$ is a cubic polynomial with roots at the vertices of a
triangle, then the roots of $p^{\prime }(z)$ are the foci of the Steiner
inellipse. A proof of this fact goes all the way back to Siebeck in 1864(see
[MP]) and is known in the literature as Marden's Theorem. Now it is easy to
show that the orthogonal least squares line, $\pounds $, from Theorem A is
identical to the line through the roots of $p^{\prime }(z)$. Thus Marden's
Theorem implies Theorem B and hence is a stronger statement than Theorem B.
There is an obvious way to try to generalize such a result to convex
quadrilaterals, D. If $p(z)$ is a quartic polynomial with roots at the
vertices of D, must the foci of the unique ellipse of maximal area inscribed
in D\ equal the roots of $p^{\prime \prime }(z)$ ? We give an example in
section 3 that shows that such a stronger statement does not hold for
parallelograms, or even for rectangles.

\section{An Area Inequality}

\textbf{Theorem 2.1:} Let $E$ be any ellipse inscribed in a convex
quadrilateral, D. Then $\dfrac{\text{Area}\left( E\right) }{\text{Area(D)}}%
\leq \dfrac{\pi }{4}$, and equality holds if and only if D\ is a
parallelogram and $E$ is tangent to the sides of D\ at the midpoints.

Before proving Theorem 2.1, we need the following lemma.

\textbf{Lemma 2.2: }Suppose that $s$ and $t$ are positive real numbers with $%
s+t>1$ and $s\neq 1\neq t$. Let

$h_{a}=\dfrac{1}{6}\dfrac{1}{t-1}\left( \allowbreak st+t-2s-1+\sqrt{%
\allowbreak (t-1)^{2}+s^{2}(t^{2}-t+1)-s(t^{2}-3t+2)}\right) $. Then $%
h_{a}\in I=$ the open interval with endpoints $\dfrac{1}{2}$ and $\dfrac{1}{2%
}s$.

\textbf{Proof of Lemma 1:}%
\begin{equation}
h_{a}-\dfrac{1}{2}=\allowbreak \dfrac{st-2\left( s+t-1\right) +\sqrt{%
(t-1)^{2}+s^{2}(t^{2}-t+1)-s(t^{2}-3t+2)}}{6\left( t-1\right) },\text{ and} 
\tag{1}
\end{equation}

\begin{equation}
h_{a}-\dfrac{1}{2}s=\allowbreak \dfrac{1}{6\left( t-1\right) }\left( \sqrt{%
(t-1)^{2}+s^{2}(t^{2}-t+1)-s(t^{2}-3t+2)}-\left( 2st-\left( s+t-1\right)
\right) \right) .  \tag{2}
\end{equation}

There are four cases to consider: $s,t>1$, $s<1<t$, $t<1<s$, and $s,t<1$. We
prove the first two cases, the proof of the other two being similar.

\textbf{Case 1:} $s,t>1$, which implies that $I=\left( \dfrac{1}{2},\dfrac{1%
}{2}s\right) $

By (1), $h_{a}-\dfrac{1}{2}>0\iff $

$\sqrt{(t-1)^{2}+s^{2}(t^{2}-t+1)-s(t^{2}-3t+2)}>2\left( s+t-1\right) -st$,
which always holds since

$\allowbreak (t-1)^{2}+s^{2}(t^{2}-t+1)-s(t^{2}-3t+2)-\left( 2\left(
s+t-1\right) -st\right) ^{2}=\allowbreak 3\left( t-1\right) \left(
s-1\right) \left( s+t-1\right) >0$

By (2), $h_{a}-\dfrac{1}{2}s<0\iff $

$\sqrt{(t-1)^{2}+s^{2}(t^{2}-t+1)-s(t^{2}-3t+2)}<2st-\left( s+t-1\right) $,
which always holds since

$\left( 2st-\left( s+t-1\right) \right) ^{2}-\left(
(t-1)^{2}+s^{2}(t^{2}-t+1)-s(t^{2}-3t+2)\right) =\allowbreak 3st\left(
t-1\right) \left( s-1\right) >0$

\textbf{Case 2:} $s<1<t$, which implies that $I=\left( \dfrac{1}{2}s,\dfrac{1%
}{2}\right) $

By (1), $h_{a}-\dfrac{1}{2}<0\iff $

$\sqrt{(t-1)^{2}+s^{2}(t^{2}-t+1)-s(t^{2}-3t+2)}<2\left( s+t-1\right) -st$,
which always holds since

$\allowbreak \left( 2\left( s+t-1\right) -st\right) ^{2}-\left(
(t-1)^{2}+s^{2}(t^{2}-t+1)-s(t^{2}-3t+2)\right) =\allowbreak 3\left(
t-1\right) \left( 1-s\right) \left( s+t-1\right) >0$

By (2), $h_{a}-\dfrac{1}{2}s>0\iff $

$\sqrt{(t-1)^{2}+s^{2}(t^{2}-t+1)-s(t^{2}-3t+2)}>2st-\left( s+t-1\right) $,
which always holds since

$\left( (t-1)^{2}+s^{2}(t^{2}-t+1)-s(t^{2}-3t+2)\right) -\left( 2st-\left(
s+t-1\right) \right) ^{2}=\allowbreak 3st\left( t-1\right) \left( 1-s\right)
>0$

\endproof%

\textbf{Proof of Theorem 2.1: }We shall prove Theorem 2.1 when D\ is not a
trapezoid, though it certainly holds in that case as well. Since ratios of
areas of ellipses and four--sided convex polygons are preserved under
one--one affine transformations, we may assume, without loss of generality,
that the vertices of D\ are $(0,0),(1,0),(0,1)$, and $(s,t)$ for some
positive real numbers $s$ and $t$. Furthermore, since D\ is convex and is
not a trapezoid, it follows easily that 
\begin{equation}
s+t>1\text{ and }s\neq 1\neq t.  \tag{3}
\end{equation}

It also follows easily that Area(D) $=\allowbreak \dfrac{1}{2}(s+t)$. Let $E$
denote any ellipse inscribed in D\ and $A_{E}=$ Area$\left( E\right) $. The
midpoints of the diagonals of D\ are $M_{1}=\left( \dfrac{1}{2},\dfrac{1}{2}%
\right) $ and $M_{2}=\left( \dfrac{1}{2}s,\dfrac{1}{2}t\right) $. Let $I$
denote the open interval with endpoints $\dfrac{1}{2}$ and $\dfrac{1}{2}s$.
By Theorem C the locus of centers inscribed in D\ is precisely the set $%
\left\{ (h,L(h),h\in I\right\} $, where 
\begin{equation}
L(x)=\dfrac{1}{2}\dfrac{s-t+2x(t-1)}{s-1}  \tag{4}
\end{equation}%
is the line thru $M_{1}$ and $M_{2}$. In proving Theorem D[H1, Theorem 3.3],
we also showed that 
\begin{equation}
A_{E}^{2}=\dfrac{\pi ^{2}}{4\left( s-1\right) ^{2}}A(h),A(h)=\allowbreak
\left( 2h-1\right) \left( s-2h\right) \left( s+2h(t-1)\right) ,h\in I. 
\tag{5}
\end{equation}%
Setting $A^{\prime }(h)=0$ yields $h=h_{a}=\dfrac{1}{6\left( t-1\right) }%
\left( \allowbreak st+t-2s-1\pm \sqrt{\allowbreak
(t-1)^{2}+s^{2}(t^{2}-t+1)-s(t^{2}-3t+2)}\right) $. By Lemma 2.2, the
positive square root always yields $h_{a}\in I$. Thus 
\begin{equation*}
h_{a}=\dfrac{1}{6\left( t-1\right) }\left( \allowbreak st+t-2s-1+\sqrt{%
\allowbreak (t-1)^{2}+s^{2}(t^{2}-t+1)-s(t^{2}-3t+2)}\right)
\end{equation*}

is the unique point in $I$ such that $A(h)\leq A\left( h_{a}\right) ,h\in I$%
. Let 
\begin{equation*}
b(s,t)=\left( st-(s+t-1)\right) ^{2}+st(s+t-1).
\end{equation*}

A simple computation using (5) shows that 
\begin{equation*}
A\left( h_{a}\right) =\dfrac{1}{27\left( t-1\right) ^{2}}\left( 2ts-s-t+1-%
\sqrt{b(s,t)}\right) \allowbreak \left( ts-2t-2s+2+\sqrt{b(s,t)}\right)
\left( s+ts+t-1+\sqrt{b(s,t)}\right) .
\end{equation*}

Thus by $\allowbreak $(5) and the fact that Area(D) $=\allowbreak \dfrac{1}{2%
}(s+t)$, we have

\begin{equation}
\dfrac{A_{E}^{2}}{(\text{Area(D)})^{2}}=\dfrac{\pi ^{2}}{27}\dfrac{\left(
2ts-s-t+1-\sqrt{b(s,t)}\right) \allowbreak \left( ts-2t-2s+2+\sqrt{b(s,t)}%
\right) \left( s+ts+t-1+\sqrt{b(s,t)}\right) }{\left( s-1\right) ^{2}\left(
t-1\right) ^{2}(s+t)^{2}}.  \tag{6}
\end{equation}

Clearly $\dfrac{A_{E}^{2}}{(\text{Area(D)})^{2}}$ is \textbf{symmetric} in $%
s $ and $t$. Thus we may assume, without loss of generality, that 
\begin{equation}
s\leq t.  \tag{7}
\end{equation}

It is convenient to make the change of variable%
\begin{equation}
u=s+t-1,v=st.  \tag{8}
\end{equation}

Solving (8) for $s$ and $t$ yields $s=u+1-\dfrac{1}{2}\left( u+1\pm \sqrt{%
(u+1)^{2}-4v}\right) ,t=\dfrac{1}{2}\left( u+1\pm \sqrt{(u+1)^{2}-4v}\right) 
$. By (7),%
\begin{equation}
s=\dfrac{1}{2}\left( u+1-\sqrt{(u+1)^{2}-4v}\right) ,t=\dfrac{1}{2}\left(
u+1+\sqrt{(u+1)^{2}-4v}\right)  \tag{9}
\end{equation}

Substituting for $s$ and $t$ using (9) gives $b(s,t)=c(u,v)$, where $%
c(u,v)=\allowbreak u^{2}-uv+v^{2}$. Then

$\dfrac{A_{E}^{2}}{(\text{Area(D)})^{2}}=\left( \dfrac{\pi ^{2}}{108\left(
v-u\right) ^{2}}(2v-u-\sqrt{\allowbreak c(u,v)})(v-2u+\sqrt{\allowbreak
c(u,v)})(v+u+\sqrt{\allowbreak c(u,v)})\right) /\left( \dfrac{1}{4}%
(u+1)^{2}\right) $, which simplifies to

\begin{equation}
\dfrac{A_{E}^{2}}{(\text{Area(D)})^{2}}=\dfrac{\pi ^{2}}{27}d(u,v)\text{,
where}  \tag{10}
\end{equation}

\begin{equation}
d(u,v)=\dfrac{\left( v-2u\right) \left( 2v-u\right) \left( u+v\right)
+2\left( c(u,v)\right) ^{3/2}}{(u+1)^{2}\left( v-u\right) ^{2}}.  \tag{11}
\end{equation}

We consider two cases: $s>1$ and $0<s<1$

\textbf{Case 1:} $s>1$

Let $w=\dfrac{u}{v}$. Then by (11) and some simplification, we have 
\begin{equation}
d(u,v)=\dfrac{v}{(u+1)^{2}}z(w),  \tag{12}
\end{equation}%
where 
\begin{equation*}
z(w)=\dfrac{\left( 1-2w\right) \left( 2-w\right) \left( 1+w\right) +2\left(
w^{2}-w+1\right) ^{3/2}}{\left( 1-w\right) ^{2}}.
\end{equation*}

Note that by (3), $u>0$, and $(u+1)^{2}-4v=(s+t)^{2}-4st=\allowbreak \left(
s-t\right) ^{2}\geq 0$, which implies 
\begin{equation}
\dfrac{v}{(u+1)^{2}}\leq \dfrac{1}{4}  \tag{13}
\end{equation}

Also, $1<s\iff 2<u+1-\sqrt{(u+1)^{2}-4v}\iff $

$\sqrt{(u+1)^{2}-4v}<u-1\iff (u+1)^{2}-4v<\left( u-1\right) ^{2}\iff u<v$,
which implies that $0<w<1$

$z^{\prime }(w)=\allowbreak \dfrac{2w^{3}-6w^{2}+9w-1+\left(
2w^{2}-5w-1\right) \sqrt{w^{2}-w+1}}{\left( w-1\right) ^{3}}\allowbreak $

$z^{\prime }(w)=0\Rightarrow \left( 2w^{3}-6w^{2}+9w-1\right) ^{2}-\left(
2w^{2}-5w-1\right) ^{2}\left( w^{2}-w+1\right) =0\Rightarrow \allowbreak
27w\left( w-1\right) ^{3}=0$, which has no solution in $(0,1)$.

$z(0)=\allowbreak 4$ and $\lim\limits_{w\rightarrow 1}z(w)=\allowbreak 
\dfrac{27}{4}$. Thus 
\begin{equation}
z(w)<\dfrac{27}{4},0<w<1,  \tag{14}
\end{equation}%
which implies, by (13), that $\dfrac{v}{(u+1)^{2}}z(w)\leq \dfrac{1}{4}%
\dfrac{27}{4}=\allowbreak \dfrac{27}{16}$. $\dfrac{\text{Area}\left(
E\right) }{\text{Area(D)}}<\dfrac{\pi }{4}$ then follows immediately by (10)
and (12).

\textbf{Case 2:} $0<s<1$

Now $s<1$ implies that $v<u$ and $v>0$ since $s,t>0$. Now we let $w=\dfrac{v%
}{u}$, which again implies that $0<w<1$. Then by (11) and some
simplification, we have 
\begin{equation}
d(u,v)=\dfrac{u}{(u+1)^{2}}z(w).  \tag{15}
\end{equation}

$\left( u-1\right) ^{2}\geq 0\Rightarrow (u+1)^{2}-4u\geq 0\Rightarrow $%
\begin{equation}
\dfrac{u}{(u+1)^{2}}\leq \dfrac{1}{4}  \tag{16}
\end{equation}

By (14) and (16), $\dfrac{u}{(u+1)^{2}}z(w)\leq \allowbreak \dfrac{27}{16}$. 
$\dfrac{\text{Area}\left( E\right) }{\text{Area(D)}}<\dfrac{\pi }{4}$ then
follows immediately by (10) and (15).

We have proven that $\dfrac{\text{Area}\left( E\right) }{\text{Area(D)}}<%
\dfrac{\pi }{4}$when D\ is not a trapezoid. Using a limiting argument, one
can then show immediately that $\dfrac{\text{Area}\left( E\right) }{\text{%
Area(D)}}\leq \dfrac{\pi }{4}$ when D\ is a trapezoid. However, that still
does not give the strict inequality when D\ is not a parallelogram. We shall
omit the details here, but the author will provide them upon request. To
finish the proof of Theorem 2.1, we need to consider the case when D\ is a
parallelogram. That case will follow from Theorem 2.3 below. 
\endproof%

\qquad First we prove that, just like any triangle in the plane, there is an
ellipse tangent to a parallelogram at the midpoints of its sides.

\textbf{Proposition 2.2:} The ellipse of maximal area inscribed in a
parallelogram, P, is tangent to P\ at the midpoints of the four sides.

\textbf{Proof:} Proposition 2.2 was proven in [H2, Theorem 3] when P is a
rectangle. The reader can prove it directly themselves using Lemma F. Ratios
of areas of ellipses, points of tangency, and midpoints of line segments are
preserved under one--one affine transformations. Since any given
parallelogram is affinely equivalent to a rectangle, that proves Proposition
2.2 for parallelograms in general. 
\endproof%

\textbf{Theorem 2.3: }Let $E$ be any ellipse inscribed in a parallelogram,
P. Then $\dfrac{\text{Area}\left( E\right) }{\text{Area(P)}}\leq \dfrac{\pi 
}{4}$, with equality if and only if $E$ is the midpoint ellipse.

\textbf{Proof: }As noted earlier, ratios of areas of ellipses and
four--sided convex polygons, are preserved under one--one affine
transformations. In addition, points of tangency and midpoints of line
segments are also preserved under one--one affine transformations. Thus we
may assume that P is a rectangle(or even a square, though we don't need that
much simplification here), $Z$, with vertices $(0,0),(l,0),(0,k),$and $(l,k),
$ where $l,k>0$. Let $E$ be any ellipse inscribed in $Z$, and let $a$ and $b$
denote the lengths of the semi--major and semi--minor axes, respectively, of 
$E$.

By Lemma F(C), 
\begin{gather*}
a^{2}b^{2}=\dfrac{2l^{2}\left( k-v\right) v}{k^{2}+l^{2}-\sqrt{%
(k^{2}+l^{2})^{2}-16l^{2}\left( k-v\right) v}}\times \\
\dfrac{2l^{2}\left( k-v\right) v}{k^{2}+l^{2}+\sqrt{%
(k^{2}+l^{2})^{2}-16l^{2}\left( k-v\right) v}}=
\end{gather*}%
\begin{gather*}
\dfrac{4l^{4}\left( k-v\right) ^{2}v^{2}}{\left( k^{2}+l^{2}\right)
^{2}-\left( (k^{2}+l^{2})^{2}-16l^{2}\left( k-v\right) v\right) }= \\
\dfrac{4l^{4}\left( k-v\right) ^{2}v^{2}}{\allowbreak 16l^{2}\left(
k-v\right) v}=\allowbreak \dfrac{1}{4}l^{2}\left( k-v\right) v.
\end{gather*}

Hence Area$\left( E\right) =\pi a^{2}b^{2}=\dfrac{\pi }{2}l\sqrt{v\left(
k-v\right) }$. Since the function $v\left( k-v\right) $ attains its unique
maximum on $(0,v)$ when $v=\dfrac{k}{2}$, the unique ellipse of maximal
area, $E_{A}$, satisfies Area$\left( E_{A}\right) =\dfrac{\pi }{2}l\sqrt{%
\dfrac{k}{2}\left( k-\dfrac{k}{2}\right) }=\dfrac{\pi }{4}lk\allowbreak $,
and thus $\dfrac{\text{Area}\left( E_{A}\right) }{\text{Area}\left( Z\right) 
}=\dfrac{\pi }{4}$. Furthermore, $\dfrac{\text{Area}\left( E\right) }{\text{%
Area}\left( Z\right) }<\dfrac{\pi }{4}$ \ when $E\neq E_{A}$. By Lemma F(B),
letting $v=\dfrac{k}{2}$ gives the points of tangency $\left( \dfrac{l}{2}%
,0\right) ,\left( 0,\dfrac{k}{2}\right) ,\left( \dfrac{l}{2},k\right) ,$and\ 
$\left( l,\dfrac{k}{2}\right) $, which are the midpoints of $Z$. That
completes the proof of Theorem 2.3 as well as the equality part of Theorem
2.1. 
\endproof%

\section{Orthogonal least squares \& Inscribed Ellipses}

If $l$ is a line in the plane, we let $d\left( z_{k},l\right) $ denote the
distance from $z_{k}$ to $l$.

\textbf{Theorem 3.1: }Let P\ be the parallelogram with vertices $%
z_{1},z_{2},z_{3},z_{4}$

(A) If P is not a square, then there is a unique line, $\pounds $, which
minimizes $\tsum\limits_{k=1}^{4}d^{2}\left( z_{k},l\right) $ among all
lines, $l$. Furthermore, the foci of the midpoint ellipse for P lie on $%
\pounds $.

(B) If P is a square, then the midpoint ellipse is a circle, and every line
through the center of P is a line of best fit for the vertices of P.

\textbf{Remarks:} (1) For part (A), an equivalent statement is that the foci
of the midpoint ellipse for P lie on the line through the second derivative
of the polynomial with roots at the vertices of D.

(2) Note that $\pounds $ is \textbf{not} the standard regression line in
general.

$\allowbreak $\textbf{Proof: }The line through the foci of an ellipse is 
\textbf{not} preserved in general under nonsingular affine transformations
of the plane, but it is preserved under translations and rotations. Thus we
can assume that the vertices of P are $O=(0,0),P=(l,0),Q=(d,k),$ and $%
R=(l+d,k)$, where $l,k>0,d\geq 0$. Using complex notation for the vertices, $%
z_{1}=0$, $z_{2}=l$, $z_{3}=d+ki$, and $z_{4}=l+d+ki$, the centroid, $g$, of
P is given by $g=\dfrac{1}{4}\tsum\limits_{k=1}^{4}z_{k}=\allowbreak \dfrac{1%
}{2}\left( l+d+ki\right) $, and $Z=\tsum\limits_{j=1}^{4}\left(
z_{j}-g\right) ^{2}=\tsum\limits_{j=1}^{4}\left( z_{j}-\dfrac{1}{2}\left(
l+d+ki\right) \right) ^{2}$, which simplifies to 
\begin{equation}
Z=d^{2}+l^{2}-k^{2}+2idk.  \tag{17}
\end{equation}

It is well known, and easy to show, that the center of any ellipse inscribed
in P equals the point of intersection of the diagonals of P, which is $%
\left( \dfrac{1}{2}\left( d+l\right) ,\dfrac{1}{2}k\right) $. Also, as noted
in the proof of Theorem 2.3, $v=\dfrac{k}{2}$ yields the unique ellipse of
maximal area, $E_{A}$, inscribed in P. Let the major axis line of $E_{A}$
refer to the line through the foci of $E_{A}$. Thus the major axis line of $%
E_{A}$ passes through $g$.

We take care of some special cases first.

\textbullet\ $d=0$: Then P is a rectangle. If $l=k$, then P is a square and
letting $v=\dfrac{k}{2}$ shows that $E_{A}$ is a circle. Also, $Z=0$. Thus
every line through the center of P is a line of best fit for the vertices of
P. We may assume now that $l\neq k$, which implies that $Z\neq 0$. Since $g=%
\dfrac{1}{2}\left( l+ki\right) $, the major axis line of $E_{A}$ passes
through $g$; $Z=l^{2}-k^{2}$; If $l>k$, then $\sqrt{Z}=\sqrt{l^{2}-k^{2}}$
is real, and thus the line, $\pounds $, thru $g$ parallel to the vector from 
$(0,0)$ to $\sqrt{Z}$ has slope $0$. By Lemma F(A), the major axis line of
any ellipse inscribed in P is parallel to the $x$ axis. Since the major axis
line passes through $g$, it must be identical to $\pounds $. If $l<k$, then $%
\sqrt{Z}=\sqrt{l^{2}-k^{2}}$ is imaginary, and thus $\pounds $ is vertical.
By Lemma F(A) again, the major axis line of any ellipse inscribed in P is
parallel to the $y$ axis. Again, the the major axis line must be identical
to $\pounds $.

\textbullet\ Assume now that $d\neq 0$, which implies that $Z\neq 0$ since $%
dk\neq 0$.

Letting $v=\dfrac{k}{2}$ in Proposition G gives the coefficients for the
equation of $E_{A}$. In particular, 
\begin{equation}
A=k^{3},B=\allowbreak k\left( d^{2}+l^{2}\right) ,C=-k^{2}d.  \tag{18}
\end{equation}%
Note that $C<0$ for all $d,k,$ and $l$, and $A<B\iff k^{3}<k\left(
d^{2}+l^{2}\right) \iff d^{2}+l^{2}-k^{2}>0$ since $k>0$.

Let $\phi $\ denote the counterclockwise angle of rotation from the line
thru the center parallel to the $x$\ axis to the major axis of $E_{A}$.

We use the formula 
\begin{equation}
\dfrac{\func{Im}\sqrt{Z}}{\func{Re}\sqrt{Z}}=\dfrac{\left\vert Z\right\vert -%
\func{Re}Z}{\func{Im}Z},\func{Im}Z\neq 0.\   \tag{19}
\end{equation}

If $k^{2}-d^{2}-l^{2}=0$, then by (17) and (19), $Z=2idk\Rightarrow \dfrac{%
\func{Im}\sqrt{Z}}{\func{Re}\sqrt{Z}}=\dfrac{2dk}{2dk}=1$. By (18), $A=B$
and $C<0$, which implies, by Lemma E, that $\phi =\dfrac{\pi }{4}$. Thus the
major axis line has slope equal to $\tan \phi =\tan \dfrac{\pi }{4}%
=\allowbreak 1$, which is the same slope as that of $\pounds $. Since both
the major axis line and $\pounds $ pass through $g$, they must be identical.

\qquad Assume now that $k^{2}-d^{2}-l^{2}\neq 0$, which implies that $A\neq
B $.

$\left\vert Z\right\vert ^{2}=\left( d^{2}+l^{2}-k^{2}\right) ^{2}+\left(
2dk\right) ^{2}=\allowbreak \left( (k+l)^{2}+d^{2}\right) \left(
(k-l)^{2}+d^{2}\right) $, which implies that

\begin{equation*}
\left\vert Z\right\vert =\sqrt{\left( (k+l)^{2}+d^{2}\right) \left(
(k-l)^{2}+d^{2}\right) }.
\end{equation*}
Hence 
\begin{equation}
\dfrac{\left\vert Z\right\vert -\func{Re}Z}{\func{Im}Z}=\dfrac{\sqrt{\left(
(k+l)^{2}+d^{2}\right) \left( (k-l)^{2}+d^{2}\right) }-\left(
d^{2}+l^{2}-k^{2}\right) }{2dk}.  \tag{20}
\end{equation}

By Lemma E, $\cot 2\phi =\dfrac{A-B}{2C}$ or $\cot (2\phi -\pi )=\dfrac{A-B}{%
2C}$, which implies that $\tan 2\phi =\dfrac{2C}{A-B}$ or $\tan (2\phi -\pi
)=\dfrac{2C}{A-B}$. Since $\tan (x-\pi )=\allowbreak \tan x$, in either
case, 
\begin{equation*}
\dfrac{2\tan \phi }{1-\tan ^{2}\phi }=\dfrac{2C}{A-B},
\end{equation*}

which implies that 
\begin{eqnarray*}
\dfrac{\tan \phi }{1-\tan ^{2}\phi } &=&\dfrac{C}{A-B}= \\
\dfrac{-k^{2}d}{k^{3}-k\left( d^{2}+l^{2}\right) } &=&\allowbreak \dfrac{dk}{%
d^{2}+l^{2}-k^{2}}
\end{eqnarray*}

by (18).

The equation$\dfrac{\tan \phi }{1-\tan ^{2}\phi }=r$ has solution 
\begin{equation}
\tan \phi =\dfrac{1}{2r}\left( -1\pm \sqrt{1+4r^{2}}\right) .  \tag{21}
\end{equation}

Letting $r=\dfrac{dk}{d^{2}+l^{2}-k^{2}}$ yields $1+4r^{2}=\dfrac{\left(
d^{2}+(k-l)^{2}\right) \left( d^{2}+(k+l)^{2}\right) }{\left(
k^{2}-d^{2}-l^{2}\right) ^{2}}$. By (21), 
\begin{equation*}
\tan \phi =\dfrac{k^{2}-d^{2}-l^{2}}{2dk}\left( 1\pm \dfrac{\sqrt{\left(
(k-l)^{2}+d^{2}\right) \left( (k+l)^{2}+d^{2}\right) }}{k^{2}-d^{2}-l^{2}}%
\right) .
\end{equation*}

There are two cases to consider.

\textbf{Case 1:} $k^{2}-d^{2}-l^{2}>0$

Let $S=\left\{ (d,k,l:d>0,k>0,l>0\right\} $, which is a connected set.

$k^{2}-d^{2}-l^{2}=0\Rightarrow \phi =\dfrac{\pi }{4}\Rightarrow \tan \phi
>0 $. Also, $\tan \phi =0\iff r=0$. But $r\neq 0$ since $dk\neq 0$. Thus $%
\tan \phi >0$ on $S$.

Also, \ 
\begin{equation*}
4k^{2}d^{2}>0
\end{equation*}

implies that

\begin{equation*}
\left( (k-l)^{2}+d^{2}\right) \left( (k+l)^{2}+d^{2}\right) -\left(
k^{2}-d^{2}-l^{2}\right) ^{2}>0,
\end{equation*}

which implies that

\begin{equation*}
\sqrt{\left( (k-l)^{2}+d^{2}\right) \left( (k+l)^{2}+d^{2}\right) }%
>k^{2}-d^{2}-l^{2}.
\end{equation*}%
Hence 
\begin{equation}
\tan \phi =\dfrac{k^{2}-d^{2}-l^{2}}{2dk}\left( 1+\dfrac{\sqrt{\left(
(k-l)^{2}+d^{2}\right) \left( (k+l)^{2}+d^{2}\right) }}{k^{2}-d^{2}-l^{2}}%
\right)  \tag{23}
\end{equation}%
since $\dfrac{k^{2}-d^{2}-l^{2}}{2dk}\left( 1-\dfrac{\sqrt{\left(
(k-l)^{2}+d^{2}\right) \left( (k+l)^{2}+d^{2}\right) }}{k^{2}-d^{2}-l^{2}}%
\right) <0.$

\textbf{Case 2:} $k^{2}-d^{2}-l^{2}<0$

Then $\dfrac{k^{2}-d^{2}-l^{2}}{2dk}\left( 1-\dfrac{\sqrt{\left(
(k-l)^{2}+d^{2}\right) \left( (k+l)^{2}+d^{2}\right) }}{k^{2}-d^{2}-l^{2}}%
\right) <0$, so again (23) holds.

By (20), 
\begin{gather*}
\dfrac{\left\vert Z\right\vert -\func{Re}Z}{\func{Im}Z}=\dfrac{%
k^{2}-d^{2}-l^{2}}{2dk}\times \\
1+\dfrac{\sqrt{\left( (k-l)^{2}+d^{2}\right) \left( (k+l)^{2}+d^{2}\right) }%
}{k^{2}-d^{2}-l^{2}}.
\end{gather*}
Thus the major axis line and $\pounds $ each have slope given by the right
hand side of (23). Since both the major axis line and $\pounds $ pass
through $g$, they must be identical. That completes the proof of (A). (B)
follows immediately from the fact that the ellipse of maximal area inscribed
in a square is a circle. 
\endproof%

We show in the following example that a version of Marden's Theorem does not
hold for ellipses inscribed in parallelograms.

\textbf{Example:} Let D\ be the rectangle with vertices $z_{1}=0$, $z_{2}=1$%
, $z_{3}=1+2i$, $z_{4}=2i$. Letting $l=1$, $k=2$, and $v=\dfrac{k}{2}=1$ in
Lemma F(A) yields the equation $4x^{2}+y^{2}-4x-2y+1=0$ for the maximal area
ellipse, $E_{A}$, inscribed in D. Rewriting the equation as $4\left( x-%
\dfrac{1}{2}\right) ^{2}+\left( y-1\right) ^{2}=1$ shows that the foci of $%
E_{A}$ are $\left( \dfrac{1}{2},1+\dfrac{\sqrt{3}}{2}\right) $ and $\left( 
\dfrac{1}{2},1-\dfrac{\sqrt{3}}{2}\right) $. If 
\begin{equation*}
Q(z)=z(z-1)(z-(1+2i))(z-2i)),
\end{equation*}%
then the roots of $Q^{\prime \prime }(z)$ are $\dfrac{1}{2}+\dfrac{1}{2}i$
and $\dfrac{1}{2}+\dfrac{3}{2}i$, which would give the points $\left( \dfrac{%
1}{2},\dfrac{1}{2}\right) $ and $\left( \dfrac{1}{2},\dfrac{3}{2}\right) $.
Hence the roots of $Q^{\prime \prime }(z)$ do \textbf{not} give the foci of
the maximal area ellipse.

\textbf{Remark: }Some heuristic reasoning shows that a version of Marden's
Theorem does not hold for ellipses inscribed in convex quadrilaterals. For
if it did, then when the quadrilateral is tangential(contains an inscribed
circle), $Q^{\prime \prime }$ would have to have a double root. But a square
is tangential, while $Q^{\prime \prime }$ has distinct roots.

For ellipses circumscribed about a convex quadrilateral, D(that is, passing
through the vertices of D), we make a conjecture similar to Theorem 2.1.

\textbf{Conjecture:} Let $E$ be any ellipse circumscribed about a convex
quadrilateral, D. Then $\dfrac{A_{E}}{\text{Area(D)}}\geq \dfrac{\pi }{2}$,
with equality if and only if D\ is a parallelogram.

There is strong evidence that this conjecture is true, but we have not yet
been able to prove this inequality.

\section{References}

[CH] G. D. Chakerian, A Distorted View of Geometry, MAA, Mathematical Plums,
Washington, DC, 1979, 130-150.

[H1] Alan Horwitz, \textquotedblleft Ellipses of maximal area and of minimal
eccentricity inscribed in a convex quadrilateral\textquotedblright ,
Australian Journal of Mathematical Analysis and Applications, 2(2005), 1-12.

[H2] Alan Horwitz, Ellipses inscribed in parallelograms, submitted to the
Australian Journal of Mathematical Analysis and Applications.

[K] D. Kalman, An elementary proof of Marden's theorem, American
Mathematical Monthly 115(2008), 330--338.

[MP] D. Minda and S. Phelps, Triangles, Ellipses, and Cubic Polynomials,
American Mathematical Monthly 115(2008), 679--689.

[W] Weisstein, Eric W. "Ellipse." From MathWorld--A Wolfram Web Resource.
http://mathworld.wolfram.com/Ellipse.html

\end{document}